\documentclass[12pt]{article}
\usepackage{amsthm,amstext, amsmath,latexsym,amsbsy,amssymb}
\usepackage{amsfonts,graphicx}
\usepackage[colorlinks=true]{hyperref}

\newtheorem{proposition}{Proposition}[section]

\numberwithin{equation}{section}
\newtheorem{theorem}{Theorem}[section] 

\def \ep{\hfill{\tiny {$\blacksquare$}}\\}
\def \bp{\noindent{\bf Proof:}\;}

\allowdisplaybreaks

\title{On the exponential convergence rate for a non-gradient Fokker-Planck equation in Computational Neuroscience}
\author{J-A. Carrillo$^1$, S. Mancini$^2$, M.-B. Tran$^3$\\
1. \it Department of Mathematics,Imperial College London, \\
\it London SW7 2AZ, U.K.\\
2. \it MAPMO, UMR 7349, CNRS, Universit\'e d'Orl\'eans, \\
\it Route de Chartres, B.P. 6759, 45067 Orl\' eans cedex 2, France\\
3. \it Department of Mathematics University of Wisconsin Madison, \\
\it WI 53706, U.S.A.}
\date{}

\begin{document}

\maketitle

\begin{abstract}
This paper concerns the proof of the exponential rate of convergence of the solution of  a Fokker-Planck equation, {with a drift term not being the gradient of a potential function and endowed by Robin type boundary conditions. This kind of problem arises, for example, in the study of interacting neurons populations.} Previous studies have numerically shown that, after a small period of time, the solution of the evolution problem  exponentially converges to the stable state of the equation. 
\end{abstract}

\section{Introduction}\label{intro}

The exponential convergence of the solutions to evolution problems towards the steady states has been largely addressed for many years, see for example \cite{AlikakosBates:SOS:1988,AMTU,CJMTU,DesvillettesVillani:OTT:2005,Haraux:CAL:1978,Hirsch:DEA:1983,GRE}. Techniques and proofs are usually based on the nature of each equation being the general entropy method developed in \cite{GRE} for linear problems and used in computational neuroscience in \cite{CCP,CCM} a powerful method for investigation of this question in the problem under consideration here.

In this work, we are concerned with the mathematical proof of the exponential rate of convergence of the solution of a non-gradient Fokker-Planck equation towards its stationary state. This problem arises, for instance, in the modeling of the evolution of the firing rates of two population of interacting neurons. {In this framework,} the partial differential equation describing the evolution of the probability distribution function $p=p(t,\nu)$, with $\nu=(\nu_1,\nu_2)\in \Omega$ and $t\geq 0$ is the Fokker-Planck equation (or forward Kolmogorov equation):
\begin{equation*}
\partial_t p + \nabla \cdot \left( Fp - {\beta^2\over 2}\nabla p\right) =0\, , \quad \textrm{in} \ [0, \infty [\times \Omega
\end{equation*}
where $\Omega$ is a bounded domain of ${\mathbb R}^2$ and the vector field $F=F(\nu)$ is defined by:
\begin{equation*}
F=\left( 
\begin{array}{c}
-\nu_1+\phi(\lambda_1+w_{11} \nu_1 + w_{12} \nu_2)\\
-\nu_2+\phi(\lambda_2+w_{21} \nu_1 + w_{22} \nu_2)\\
\end{array}
\right)
\end{equation*}
with $\phi(z)$ a sigmoid function and $w_{ij}$ positive weights. The Fokker-Planck equation is endowed by the following Robin or no flux boundary conditions:
\begin{equation*}
\left(Fp - {\beta^2\over 2}\nabla p\right)\cdot n =0  \, , \ \textrm{on} \ \partial \Omega
\end{equation*}
with $n$ the outward normal unit vector on $\partial \Omega$, $\beta$ being the noise level, and a non-negative initial data $p(0,\nu)=p_0(\nu) \geq 0$. It is easily seen that the $F$ does not satisfy to the Schwartz conditions, that is $\partial_{\nu_2} F_1 \neq \partial_{\nu_1} F_2$, so that  there exists no potential function $V(\nu)$ such that $F=-\nabla V$. This in particular implies that it is less likely there is an explicit formulation of the associated steady state, as it would be the case with a drift $F$ defined by $F=-\nabla V$ for some smooth potential function $V$. We refer to \cite{ArnoldCarrilloManzini:RLT:2010,ArnoldCarlen:AGB:1999,ArnoldCarlenJu:LTB:2008} for some discussion on these issues and some cases in which, despite of the non-gradient drifts in the Fokker-Planck equations, one is able to find explicit stationary states. In our case, for general domains and so complicated drift it is in general not possible to find explicitly these stationary states. Note that as soon as we find a potential function $V$ such that the drift $F=-\nabla V$, the steady state reads
$$
p_\infty=C\exp\left( \frac{-2 V}{\beta^2}\right),
$$
with $C$ a normalization coefficient.

In this work, we will deal with estimates of the exponential convergence rates for general linear Fokker-Planck equations with non-gradients drifts of the form:
\begin{equation}\label{FPNeuroscience3}
\left \{
\begin{array}{ll}\partial_tu-\Delta{u}+\nabla \cdot (Fu)=0 \quad\mbox{ in } {\Omega}\times(0,{T})\vspace{.1in},\\
u(0,\cdot)=u_0(\cdot) \quad \mbox{ in } {\Omega},\vspace{.1in}\\
(Fu-\nabla u)\cdot n=0 \quad \mbox{ on } \partial{\Omega}.\end{array}\right. 
\end{equation}
where the unknown $u(t,x)$ is a probability density function, the drift $F\in C^2(\overline{\Omega})$ satisfies the incoming boundary condition 
\begin{equation}\label{hypoF}
F \cdot n <0 \quad \textrm{on}\ \partial \Omega,
\end{equation} and with the initial data normalized to 1, $\int_\Omega u_0 (x)\ dx =1$. We have assumed unit diffusion constant for simplicity without loss of generality.
Under these hypothesis, existence, uniqueness and positivity of the solution $u=u(t,x)$ of the evolution problem \eqref{FPNeuroscience3}, and of the stable state (stationary solution)  $u_\infty(x)$ of the associated problem were proved in \cite[Theorem 2]{CCM}, as well as the mass density conservation, 
$$
\int_\Omega u(t,x)\ dx =\int_\Omega u_0(x)\ dx =1,
$$
and the $L^2$-convergence of the time dependent solution to the stationary solution $u_\infty$. Numerical simulations were also performed underlying the exponential rate of convergence of the solution of the evolution problem towards the stable state. However, the theoretical proof of this exponential rate convergence was not discussed in \cite{CCM}. 
In the next section we will discuss the exponential convergence towards the steady state $u_\infty$ by direct Poincar\'e inequalities. Section 3 is devoted to an alternative approach implying relations between the sharp exponential convergence rates with the best constants for different Poincar\'e type inequalities.


\section{Convergence to equilibrium}\label{proof}
Consider a drift $F\in C^2(\overline{\Omega})$ such that $F\cdot n <0 $ and let $u=u(t,x)$ be the unique solution of problem \eqref{FPNeuroscience3}, and $u_\infty=u_\infty(x)$ be the solution of the stationary associated problem:
\begin{eqnarray}\label{stationary-u}
\begin{cases}-\Delta u_\infty+\nabla(F u_\infty)=0 \quad \mbox{ in } {\Omega}\vspace{.1in},\\
(F u_\infty -\nabla u_\infty )\cdot n = 0\quad \mbox{ on } \partial{\Omega}\,,\end{cases} 
\end{eqnarray}
ensured by \cite[Theorem 2]{CCM} and satisfying
$$
\int_\Omega  u_0(x) \, dx =\int_\Omega u_\infty(x)\,dx =1.
$$
We shall prove the following result:
\begin{theorem}\label{u-conv}
The solution $u$ to \eqref{FPNeuroscience3} exponentially converges to the steady state $u_\infty$: there exist $\alpha$, $C>0$ such that
$$
\| u -u_\infty\|_{L^2(\Omega)} \leq C \exp(-\alpha t).
$$
\end{theorem}

Let us first remind the reader that applying the General Entropy Method \cite{GRE} adapted to this problem in \cite[Theorem 4]{CCM} with $v=1$, $u_1=u$, $u_2=u_\infty$ and $H(u_1|u_2)=(u_1-u_2)^2/u_2^2$, then we have
\begin{equation}\label{entropy-decay}
\frac{d}{dt} \mathcal{H}_v(u_1|u_2) =- \mathcal{D}_v(u_1|u_2) \leq 0,
\end{equation}
where 
$$
\mathcal{H}_v(u_1|u_2) = \int_\Omega \frac{(u-u_\infty)^2}{u_\infty} \, dx
$$
and 
$$
\mathcal{D}_v(u_2|u_1) = \int_\Omega u_\infty \left| \nabla \frac{u}{u_\infty}\right|^2\, dx\,.
$$
Therefore, in other to get the exponential convergence, the key point to prove  is the following  Poincar\'e inequality:
\begin{equation}\label{poinc-V}
 \mathcal{H}_v(u|u_\infty)\!=\!\!\int_\Omega \frac{(u-u_\infty)^2}{u_\infty} \, dx \leq \mathcal{P}_\Omega (u_\infty)\!\! \int_\Omega \! u_\infty  \left| \nabla \frac{u}{u_\infty}\right|^2\!\! dx= \mathcal{P}_\Omega(u_\infty) \,\mathcal{D}_v(u|u_\infty),
\end{equation}
where $\mathcal{P}_\Omega(u_\infty)>0$ is the best constant for the Poincar\'e inequality \eqref{poinc-V}, with weight $u_\infty$ in the domain $\Omega$. In the sequel, we will not put the subindex $\Omega$ in the Poincar\'e constants for notational simplicity. In fact, from \eqref{entropy-decay} and \eqref{poinc-V} we then deduce that
$$
\frac{d}{dt} \mathcal{H}_v(u|u_\infty) + \mathcal{P}(u_\infty) \ \mathcal{H}_v(u|u_\infty) \leq 0,
$$
and applying Gronwall Lemma
$$
\mathcal{H}_v (u|u_\infty) \leq C \exp (\alpha t),
$$
with $C=\mathcal{H}_v (u_0|u_\infty)$ and $\alpha = \mathcal{P}(u_\infty)$, leading to the exponential convergence of $u$ to $u_\infty$. Notice that \eqref{poinc-V} is equivalent to
\begin{equation}\label{poinc-Vb}
\int_\Omega u_\infty \left(\frac{u}{u_\infty}-1\right)^2\, dx \leq \mathcal{P}(u_\infty) \int_\Omega u_\infty  \left| \nabla \left(\frac{u}{u_\infty}-1\right)\right|^2\, dx\,,
\end{equation}
Put $w=\frac{u}{u_\infty}-1$, then \eqref{poinc-Vb} reads
\begin{equation*}
\int_\Omega u_\infty w^2 \, dx \leq \mathcal{P}(u_\infty) \int_\Omega u_\infty  \left| \nabla w\right|^2\, dx\, .
\end{equation*}
Note that, since both $u$ and $u_\infty$ are normalized to 1, we also have the constraint:
\begin{equation*}
\int_\Omega u_\infty w\, dx=\int_\Omega u_\infty \left(\frac{u}{u_\infty}-1\right)\, dx=\int_\Omega  \left({u}-{u_\infty}\right)\, dx=0.
\end{equation*}
Inequality \eqref{poinc-V} then follows from the following result:

\begin{proposition}\label{Poincare}
Under the assumption  
$$
\int_\Omega\Phi H\, dx=0,
$$
where $H$ is bounded from above and below by positive constants,
there exists $\mathcal{P}>0$ such that the following Poincar\'e inequality holds
$$
\int_\Omega H |\Phi|^2\, dx \leq \mathcal{P}\int_\Omega H |\nabla\Phi|^2\, dx.
$$
\end{proposition}
\bp
We follow the classical proof of Poincar\'e inequality, by contradiction. 
Suppose that there exists a sequence $\{\Phi_m\}$, such that
$$
\max(H)\!\! \int_\Omega |\Phi_m|^2\, dx\geq \!\int_\Omega H |\Phi_m|^2\, dx\geq \! m\int_\Omega H |\nabla\Phi_m|^2 \, dx\geq m \min(H) \!\!\int_\Omega |\nabla\Phi_m|^2\, dx.
$$
By normalization, we can suppose that $\|\Phi_m\|_{L^2}=1$. Therefore the sequence $\{\Phi_m\}$ is bounded in $W^{1,2}(\Omega)$. By mean of the Rellich-Kondrachov Theorem, there exists a subsequence $\{\Phi_{m_j}\}$ and a function $\bar\Phi$ in $L^2(\Omega)$ such that $\{\Phi_{m_j}\}$ converges strongly to $\bar\Phi$ in $L^2(\Omega)$. Passing to the limit $m\to\infty$, we get that $\nabla\bar\Phi=0$ a.e. and $\|\bar\Phi\|_{L^2}=1$, which implies  $\bar\Phi$ is a constant. But, since $H$ is strictly positive, this contradicts with the hypothesis
$\int_\Omega\bar\Phi H\, dx =0$, concluding the proof. \ep

Although we can obviously relate the best constant $\mathcal{P}(u_\infty)$ to the classical Poincar\'e inequality best constant $\mathcal{P}(1)$ by $\mathcal{P}(u_\infty)\leq \tfrac{\max(u_\infty)}{\min(u_\infty)}\mathcal{P}(1)$, we have no information on the value of the constant obtained by this proof. In the next section, we propose an alternative proof in which we show that the constants $C$ and $\alpha$ in Theorem \ref{u-conv} are linked to the bounds of the solution $K$ of an auxiliary problem.


\section{Alternative Proof}
In this section, we propose an alternative proof of Theorem \ref{u-conv} which give us another characterization of the constants $C$ and $\alpha$. Note first that to prove the exponential convergence of $u$ to $u_\infty$, is equivalent to prove the exponential convergence of $\Phi=u/ u_\infty$ to $1$. As we will see, the proof is based on a new conservation property, see Proposition \ref{phiK-cons} below. 

Before dealing with the exponential convergence problem for the function $\Phi$, we need to define for which problem $\Phi$ is a solution. This is done in the following. 
Let us first define $\Phi$ and the initial data $\Phi_0$ as follows, since $u_\infty >0$, we define
$$
\Phi(t,x)=\frac{u(t,x)}{u_\infty(x)} ,\quad \Phi_0(x)=\frac{u_0(x)}{u_\infty(x)}\,,
$$
and let us introduce a modified drift $F^*$ as
\begin{equation}\label{Fstar}
F^*=F-2\frac{\nabla u_\infty}{u_\infty}.
\end{equation}
We note that on the boundary $\partial \Omega$, we obtain
$$
F^* \cdot n =\left(F -2 \frac{\nabla u_\infty}{u_\infty}\right) \cdot n =-F \cdot n >0 ,
$$
since $F \cdot n =(\nabla u_\infty / u_\infty) \cdot n $ and $F$ satisfies \eqref{hypoF}.

\begin{proposition}\label{Phi-problem}
The function $\Phi=u/u_\infty$ satisfies
\begin{equation}\label{FPNeuroscience5}
\left \{
\begin{array}{ll}\partial_t\Phi-\Delta \Phi+F^*\cdot\nabla\Phi=0 \mbox{ in } {\Omega}\times(0,{T})\vspace{.1in},\\
\Phi(.,0)=\Phi_0(.) \textrm{ in } {\Omega},\vspace{.1in}\\
\nabla \Phi \cdot n =0 \textrm{ on } \partial{\Omega}.\end{array}\right. 
\end{equation}
\end{proposition}

\noindent
\bp
In order to simplify notations, consider $g=1/u_\infty$,  $g\in C^2(\Omega)\cap C(\bar\Omega)$ and $g>0$, and $\Phi=u\ g$ and $\Phi_0=u_0\ g$. We derive $\Phi=u\ g$ to get, for all $i=1,2$
$$
\partial_{x_i}\Phi=g\partial_{x_i}u +u\partial_{x_i}g,
$$
$$
\partial_{x_i}^2\Phi=\partial_{x_i}^2u g+\Phi\frac{\partial_{x_i}^2g}{g}+2\frac{\partial_{x_i}g}{g}\partial_{x_i}\Phi-2\Phi\left|\frac{\partial_{x_i}g}{g}\right|^2.$$
These identities imply, for the boundary conditions on $\partial \Omega$
$$
\nabla\Phi \cdot n = \left( {\nabla u \over u_\infty} - {u\nabla u_\infty \over u_\infty^2}\right)\cdot n = \left( {\nabla u \over u_\infty} -  {uF\over u_\infty}\right)\cdot n =0,
$$
and inside the domain $\Omega$, we obtain
\begin{align}\label{eq1-bis}
\nonumber
\partial_t\Phi-\Delta\Phi&=(\partial_tu-\Delta u)g-\Phi\frac{\Delta g}{g}-2\nabla\Phi \cdot \frac{\nabla g}{g}+2\Phi\left|\frac{\nabla g}{g}\right|^2\\
\nonumber
&=-\nabla\cdot(Fu)g-\Phi\frac{\Delta g}{g}-2\nabla\Phi \cdot \frac{\nabla g}{g}+2\Phi\left|\frac{\nabla g}{g}\right|^2\\
\nonumber
&=-\nabla\cdot(F\Phi)-2\nabla\Phi \cdot \frac{\nabla g}{g}\\
&\quad +\Phi\left[-\left(\frac{\Delta g}{g}-2\left|\frac{\nabla g}{g}\right|^2\right)+F\cdot\frac{\nabla g}{g}\right].
\end{align}
Since $u_\infty=1/g,$ we have
$$
\nabla{u_\infty}=-\frac{\nabla g}{g^2} \quad \mbox{ and }\quad\Delta u_\infty=-\frac{\Delta g}{g^2}+2\frac{|\nabla g|^2}{g^3},
$$
so that, dividing \eqref{eq1-bis} by $g=1/u_\infty$, we get
\begin{align*}
u_\infty \partial_t\Phi- u_\infty \Delta \Phi&\,+ u_\infty F\cdot \nabla\Phi-2\nabla\Phi\cdot \nabla u_\infty\\
&+\Phi[-\Delta u_\infty+F\cdot \nabla u_\infty+u_\infty\nabla \cdot F ]=0.
\end{align*}
Hence $\Phi$ must satisfy
\begin{equation}\label{FPNeuroscience4}
\begin{cases}
u_\infty \partial_t\Phi-u_\infty \Delta \Phi+u_\infty F\cdot \nabla\Phi-2\nabla\Phi\cdot \nabla u_\infty \\
\qquad \qquad +\Phi[-\Delta u_\infty +\nabla\cdot (F u_\infty )]=0 \quad\mbox{ in }\ {\Omega}\times(0,{T}),\\
\Phi(\cdot,0)=\Phi_0(\cdot)\quad \mbox{ in }\ {\Omega},\\
\nabla \Phi\cdot n=0\quad \mbox{ on }\ \partial{\Omega}.
\end{cases} 
\end{equation}
Finally, recalling the definition of $F^*$ given by \eqref{Fstar}, and that $u_\infty$ solves  \eqref{stationary-u}, problem $(\ref{FPNeuroscience4})$ is converted into \eqref{FPNeuroscience5}, concluding the proof.
\ep
\\

Consider now the stationary problem associated to \eqref{FPNeuroscience5}:
\begin{equation}\label{equilibrium2}
\begin{cases}
-\Delta \Phi_\infty+F^*\cdot \nabla\Phi_\infty=0 \quad \textrm{ in }\ {\Omega}\times(0,{T})\vspace{.1in},\\
\nabla \Phi_\infty\cdot n=0 \quad \textrm{ on }\ \partial{\Omega}.
\end{cases} 
\end{equation}
then it is easily seen that constant are solutions to \eqref{equilibrium2}, so that proving the exponential convergence of $\Phi$ to $\Phi_\infty =1$, will give us the wanted exponential convergence of $u$ to $u_\infty$. So, the main results we have to prove, yielding to the exponential convergence of $u$ to $u_\infty$  is the following:
\begin{theorem}\label{Phi-conv}
The solution $\Phi$ of $\eqref{FPNeuroscience5}$ exponentially converges in time to its equilibrium state $\Phi_\infty$, solution to \eqref{equilibrium2}.
\end{theorem}
\noindent
In the sequel, in order to prove Theorem \ref{Phi-conv}, we need an auxiliary problem, which is the dual problem to \eqref{equilibrium2} given by
\begin{eqnarray}\label{pb-K}
\begin{cases}\Delta K+\nabla\cdot (K F^* )=0 \quad \mbox{ in }\ {\Omega}\vspace{.1in},\\
(\nabla K + K F^*)\cdot n=0 \quad \mbox{ on }\ \partial{\Omega},
\end{cases} 
\end{eqnarray}
for which we next prove a conservation property of the function $\Phi K$. This conservation results is a key point in order to prove a Poincar\'e inequality on the function $\Phi$ and to get in an alternative way the exponential convergence.

\begin{proposition}\label{phiK-cons}
Let $K$ be a solution of \eqref{pb-K} and $\Phi$ be a solution of 
\eqref{FPNeuroscience5}, then
$$
\int_\Omega\Phi_0 K\,dx=\int_\Omega\Phi K\,dx.
$$
\end{proposition}
\bp
Note that $-F^*\cdot n<0$ on $\partial\Omega$, thus by using the same arguments as in \cite{CCM}, this problem admits a unique strictly positive and bounded solution $K\in H^2(\Omega)$. Furthermore, we have
$$
-\Delta \Phi+F^*\nabla\Phi=-\frac{1}{K}\nabla\cdot (K\nabla \Phi)+\left( \frac{\nabla{K}}{K}+F^*\right) \nabla\Phi,
$$
so that problem $(\ref{FPNeuroscience5})$ writes as
\begin{equation}\label{FPNeuroscience6}
\left \{
\begin{array}{ll}
\displaystyle \partial_t\Phi-\frac{1}{K}\nabla \cdot (K\nabla \Phi)+\left( \frac{\nabla{K}}{K}+F^*\right)\cdot\nabla\Phi=0 \quad \mbox{ in } \ {\Omega}\times(0,{T})\vspace{.1in},\\
\Phi(.,0)=\Phi_0(.) \quad \mbox{ in } \ {\Omega},\vspace{.1in}\\
\nabla \Phi \cdot n=0 \quad \mbox{ on }\ \partial{\Omega}.\end{array}\right. 
\end{equation}
Use $K$ as a test function for $(\ref{FPNeuroscience6})$, then
\begin{align*}
0&=\partial_t\int_\Omega\Phi K\,dx+\int_\Omega({\nabla{K}}+{K}F^*)\cdot \nabla\Phi\,dx \ {-\int_{\partial \Omega} K\nabla \Phi\cdot n\,d\sigma(x)}\\
&=\partial_t\int_\Omega\Phi K\,dx-\int_\Omega\nabla\cdot({\nabla{K}}+{K}F^*)\Phi\,dx+\int_{\partial\Omega}\Phi ({\nabla{K}}+{K}F^*) \cdot n\,d\sigma(x)\\
&=\partial_t\int_\Omega\Phi K\,dx-\int_\Omega({\Delta{K}}+\nabla\cdot ({{K}F^*}))\Phi\,dx =\partial_t\int_\Omega\Phi K\,dx,
\end{align*}
which concludes the proof.\ep
\\

Since, problem \eqref{pb-K} admits a unique solution, then constants are the unique solutions, in distributional sense, of the stationary problem \eqref{equilibrium2}. Hence, up to a normalization constant, we shall consider in the sequel that $\Phi_\infty=1$. Moreover,  if we consider  $\Psi=\Phi-\Phi_\infty$, then $\Psi$ still satisfies \eqref{Phi-problem} and converges to 0, in other words its equilibrium is the null function. Therefore, renaming $\Psi=\Phi$ and in order to simplify the proof of Theorem \ref{Phi-conv}, we are reduced to prove the exponential convergence of $\Phi$ solution of \eqref{Phi-problem} to 0, which is done in Theorem \ref{main-th} below. We impose the following normalization
$$
\int_\Omega\Phi_0 K\, dx=0,
$$
so that from the conservation property in Proposition \ref{phiK-cons} we get
$$
\int_\Omega\Phi K\, dx=0.
$$
Thanks to this normalization, the proof of Theorem \ref{Phi-conv} is reduced to prove the following result.

\begin{theorem}\label{main-th}
Assuming that 
$$
\int_{\Omega}\Phi_0 K\, dx=0,
$$
then, the solution $\Phi$ of $(\ref{FPNeuroscience5})$ exponentially decays in time  towards $0$, that is,
$$
\|\Phi \|_{L^2(\Omega)}\leq {\tilde C}\exp(-\tilde{\alpha} t),
$$
where $\tilde C$ and $\tilde{\alpha}$ are estimated by
$$
{\tilde C}=\frac{1}{\min(K)} \int_\Omega |\Phi_0|^2 K\, dx\ \quad\mbox{and}\quad \tilde{\alpha}=\mathcal{P}(K)\,.
$$
\end{theorem}

\noindent
{\bf Proof:} 
Use $K\Phi$ as a test function for $(\ref{FPNeuroscience6})$. Since $K$ is bounded from above and from below by some positive constants, and since from the normalization hypothesis, the Poincar\'e inequality of Proposition \ref{Poincare} holds, then:
\begin{eqnarray}\label{phiK2}
\nonumber
0&=&\frac{1}{2}\partial_t\int_\Omega|\Phi|^2K\, dx+\int_\Omega|\nabla\Phi|^2K\, dx+\int_\Omega\Phi ({\nabla{K}}+{K}F^*)\cdot \nabla\Phi\, dx \\
\nonumber
&=&\frac{1}{2}\partial_t\int_\Omega|\Phi|^2K\, dx +\int_\Omega|\nabla\Phi|^2K\, dx +\frac{1}{2}\int_\Omega({\nabla{K}}+{K}F^*)\cdot \nabla|\Phi|^2 \, dx \\
\nonumber
&=&\frac{1}{2}\partial_t\int_\Omega|\Phi|^2K\, dx +\int_\Omega|\nabla\Phi|^2K\, dx -\frac{1}{2}\int_\Omega({\Delta{K}}+\nabla\cdot({{K}F^*}))|\Phi|^2 \, dx \\
\nonumber
&=&\frac{1}{2}\partial_t\int_\Omega|\Phi|^2K\, dx +\int_\Omega|\nabla\Phi|^2K\, dx \\
\nonumber
&\geq & \frac{1}{2}\partial_t\int_\Omega|\Phi|^2K\, dx +\mathcal{P}(K)\int_\Omega|\Phi|^2K\, dx \,.
\end{eqnarray}
Inequality \eqref{phiK2} implies:
$$
\int_\Omega|\Phi|^2K\, dx\leq \left(\int_\Omega|\Phi_0|^2K\, dx\right) \exp(-2\tilde{\alpha}t),
$$
with $\tilde{\alpha}=\mathcal{P}(K)$. With some more computations, we finally obtain
$$
\int_\Omega|\Phi|^2\, dx \leq \frac{1}{\min(K)} \int_\Omega|\Phi|^2 K\, dx \leq \left(\frac{1}{\min(K)} \int_\Omega |\Phi_0|^2 K\, dx\right) \exp(-2\tilde{\alpha}t),
$$
concluding the proof of theorem \ref{main-th}.
\ep

Finally, to conclude the proof of the time exponential convergence of $u$, solution to \eqref{FPNeuroscience3}, to the steady state $u_\infty$ given by \eqref{stationary-u}, we just have to note that, since $g=1/u_\infty$ is strictly positive and bounded, Theorem \ref{Phi-conv} leads to:
$$
\|u-u_\infty\|_{L^2(\Omega)}\leq \max(u_\infty) \tilde{C}\exp(-\tilde{\alpha} t)\,,
$$
concluding the proof of Theorem \ref{u-conv}. Notice that the decay constants might be different to the ones in Theorem \ref{u-conv} obtained in Section 2 due to the different stationary problems with solutions $K$ and $u_\infty$ used. 

\section*{Acknowledgments}
\small{JAC was partially supported by the Royal Society by a Wolfson Research Merit Award and the EPSRC grant EP/K008404/1. SM was partially supported by the KIBORD project (ANR-13-BS01-0004) funded by the French Ministry of Research. MBT was partially supported by the NSF Grant RNMS (Ki-Net) 1107444 and the ERC Advanced Grant FP7-246775 NUMERIWAVES. }


\end{document}